\documentclass[11pt]{amsart}

\usepackage[margin=1in]{geometry}

\usepackage[T1]{fontenc}
\usepackage[latin1]{inputenc}

\usepackage{amsthm, amsmath, amssymb,booktabs,xcolor,graphicx}
\usepackage{mathtools}
\usepackage[shortlabels]{enumitem}
\usepackage{scalerel,stackengine}
\usepackage[textsize=scriptsize,colorinlistoftodos]{todonotes}
\usepackage[hidelinks]{hyperref}
\usepackage[noabbrev,capitalize,nameinlink]{cleveref}
\crefformat{equation}{(#2#1#3)}
\crefrangeformat{equation}{(#3#1#4) to~(#5#2#6)}
\crefmultiformat{equation}{(#2#1#3)}
{ and~(#2#1#3)}{, (#2#1#3)}{ and~(#2#1#3)}
\usepackage{caption}
\usepackage{subcaption}

\usepackage[color,final]{showkeys}

\colorlet{refkey}{orange!50}
\colorlet{labelkey}{blue!50}

\newtheorem{theorem}{Theorem}[section]
\newtheorem{thm}[theorem]{Theorem}%[section]

\newtheorem{lemma}[theorem]{Lemma}
\newtheorem{claim}[theorem]{Claim}
\newtheorem{cor}[theorem]{Corollary}

\newtheorem*{question*}{Question}

\theoremstyle{definition}
\newtheorem{defn}[theorem]{Definition}

\newtheorem{model}{Model}

\newtheorem*{definition*}{Definition}
\newtheorem{example}[theorem]{Example}
\newtheorem*{example*}{Example}

\theoremstyle{remark}
\newtheorem{rem}[theorem]{Remark}

\numberwithin{equation}{section}

\newcommand{\ignore}[1]{}

\newcommand{\ind}{\mathbf{1}}
\renewcommand{\P}{\mathbf{P}}

\stackMath
\newcommand\widecheck[1]{
\savestack{\tmpbox}{\stretchto{
  \scaleto{
    \scalerel*[\widthof{\ensuremath{#1}}]{\kern-.6pt\bigwedge\kern-.6pt}
    {\rule[-\textheight/2]{1ex}{\textheight}}
  }{\textheight}
}{0.5ex}}
\stackon[1pt]{#1}{\scalebox{-1}{\tmpbox}}
}
\DeclareMathOperator{\ab}{abs}

\DeclareMathOperator{\rand}{\overset{\text{R}}\leftarrow}

\newcommand{\mfp}{\mathfrak{p}}

\newcommand{\EE}{\mathbb{E}}

\newcommand{\RR}{\mathbb{R}}

\newcommand{\NN}{\mathbb{N}}
\newcommand{\ZZ}{\mathbb{Z}}

\makeatletter
\newcommand\thankssymb[1]{\textsuperscript{\@fnsymbol{#1}}}
\makeatother

\title{When will (game) wars end?}

\author{Manan Bhatia, Byron Chin, Nitya Mani, Elchanan Mossel}

\address{Department of Mathematics, Massachusetts Institute of Technology, Cambridge, MA 02139, USA}
\email{\{mananb,byronc,nmani,elmos\}@mit.edu}

\begin{document}

\begin{abstract}
We study several variants of the classical card game \textit{war}. As anyone who played this game knows, the game can take some time to terminate, but it usually does. Here, we analyze a number of asymptotic variants of the game, where the number of cards is $n$, and show that all have expected termination time of order $n^2$. This is the same expected termination time as in the game where at each turn a fair coin toss decides which player wins a card, known as Gambler's Ruin and studied by Pascal, Fermat and others in the seventeenth century.
\end{abstract}

\maketitle

\section{Introduction}
The \textit{game of war} is a popular children's card game, typically played with a standard deck of cards. 
We will concentrate on the simplest setting with two players. At the beginning, a deck of cards is shuffled and split into two hands, one for each player. In a series of rounds, each player flips over their top card, and the player with the larger top card collects both cards, returning them to the bottom of his/her hand. 
Ties can be broken in a variety of ways; one common mechanism is a \textit{war} round.
Here is the description from Wikipedia~\cite{wikiwar:23}: 
\begin{quote}
If the two cards played are of equal value, then there is a "war". Both players place the next card from their pile face down and then another card face-up. The owner of the higher face-up card wins the war and adds all the cards on the table to the bottom of their hand. If the face-up cards are again equal then the battle repeats with another set of face-down/up cards. This repeats until one player's face-up card is higher than their opponent's. Most descriptions of War are unclear about what happens if a player runs out of cards during a war. In some variants, that player immediately loses. In others, the player may play the last card in their hand as their face-up card for the remainder of the war or replay the game from the beginning.
\end{quote}
Different authors of the paper have played other variants where at a war round two/three cards are placed faced down before the deciding card is drawn.  

For those who have played the game, it can feel like an exercise in patience as cards move back and forth with players' hands growing in size and shrinking dramatically many times over the course of the round. A priori, it is not even obvious that the game is expected to end in a finite amount of time; however, finiteness was established in~\cite{LR12}.
While the finiteness of the game is reassuring, it begs the following natural question: how long does it take for a game of war to end? Of course, we can simulate the game, and this will give us a good idea of the distribution of the length of the game. However, we are interested in the asymptotic question: if there are $n$ cards, how long would it take the game to terminate?

When playing the game, and ignoring wars, it may feel like that at each step, each player gains or loses a card with $1/2$ probability, and therefore the process should behave like a simple random walk. Simple random walks and the Gambler's ruin have been studied since the 17th century (see e.g. \cite{edwards1982pascal, edwards1983pascal} for an account of Pascal and Fermat's communication on the subject). This analogy suggests an expected termination time of order $n^2$. There are many ways which this can be shown, see e.g. \cite[Chapter 6]{van2016probability} for a simple recursive approach, or \cite[Section 4.8]{durrett2019probability} for an approach using the theory of martingales.

However, the war game is of course not a simple random walk. 
Consider the simplest version of this question where war is played with a single-suited deck of cards $\{1, 2, \ldots, n\}$. 
In this case, and unlike the random walk case, it is clear that the player who has $n$ (the largest card) should eventually win (if any player wins), since they cannot lose card $n$ to the other player. As a result, the dependency on $n$ of the expected time to win is not clear. 
In this work, we study several variants of war including this simple variant, and show that they all have expected termination time of order $n^2$, where $n$ is the number of cards.

\subsection{Models and Main Results}
 In~\cref{s:rand}, we consider a variant of the game where each player picks a uniformly chosen card from her hand and the two cards play against each other. In our main result, we prove that this variant has expected termination time of $O(n^2)$ for a number of rules of determining which is the winning card, including rules that take into account the remaining cards in the hand or deck. Thus, for a very general setting, the expected termination time of the game is $O(n^2)$. 

In~\cref{s:gladiator}, we attempt to bound the termination time for games where cards are drawn from the top and returned to the bottom of the winning hand at each round. We can only analyze this variant in cases where the chance of winning is proportional to the {\em value of each card}. These kinds of comparisons are known as the Bradley-Terry Model, introduced by Zermelo in ~\cite{zermelo1929berechnung} in the 1920's to study pairwise interactions. This settings also generalizes the Gladiator Puzzle discussed by Winkler \cite{winkler2003mathematical}, in which the value of a card is simply its number.

\subsection{Related Work}
Lakshtanov and Roschchina~\cite{LR12} showed that the expected time for the game of war to end with the deck of cards $\{1, 2, \ldots, n\}$ is finite, as long as there is some randomness in the order in which cards are returned to the bottom of players' hands. Further, they exhibit several infinite cycles that show that allowing players to strategically return cards to the bottom of their hands can result in infinite loops (see Figures 1 and 2 in~\cite{LR12}).

In~\cite{KP22}, some combinatorial aspects of the single-suited variant of war were studied, including how many times a player cycles through their hand over the course of the game and certain special cases in which the time for one player to win can be conclusively determined. Various other aspects of the card game war and several other variants have been studied, such as in~\cite{AT12,BR17,BK02,LA13,SPI10,WR01}.

\subsection*{Notation}
Let $[n] = \{1, 2, \ldots, n\}$ and let $\mathbf{0}, \ind$ be the all-zeros and all-ones vectors (where dimension is clear from context). We use $a \rand A$ to denote a uniformly random sample $a$ from set $A$. We let $\NN$ denote the non-negative integers and $\ZZ_{>0}$ the positive integers.

\section{Random Game}\label{s:rand}

We will consider a fairly general instantiation of war that we denote $\mfp$-war. In this game, 
there are two players denoted $A$ and $B$. Given a deck $D$ of $2n$ cards, we write $A$
for player $A$'s hand and $B$ for player $B$'s hand, i.e., $D \setminus A$. In each round, player $A$ plays a random card $a$ and player $B$ plays a random card $b$ from their respective hands. Each round's outcome is determined using a function $p_{a, b}(S)$,
where $S = A \setminus \{a\}$. 
This function determines the probability that player $A$'s card beats player $B$'s card given the two hands, and that the players drew $a$ and $b$ respectively. More concretely, player $A$ wins the two cards with probability $p_{a, b}(S)$, and otherwise, the two cards are won by player $B$. The game remains in play until all the cards are in one player's hand, and this player is declared to be the winner of the game. Note that $p_{a, b}(S)$ may depend on $a$ and $b$,
$A \setminus \{a\}$ and $B \setminus \{b\}$,
which are the cards in play, as well at the remaining cards in player A's hand and player B's hand respectively. We desire for players A and B to be 
exchangeable in the game. As a result, $p$ must satisfy the condition $p_{b, a}(D \setminus (S \cup \{a, b\})) = 1 - p_{a, b}(S)$.  

\begin{defn}
    Consider a multi-set $D$ of size $2n$ with each entry being an integer. 
    A \textit{winning rule} $\mfp$ is function 
    that maps $a \neq b \in D$ and 
    $S \subseteq D \setminus \{a,b\}$ to $p_{a,b}(S)$, where $p_{a,b}(S) \geq 0$ and 
    $p_{a,b}(S) + p_{b,a}(D \setminus (S \cup \{a,b\} ) = 1$. 
    The winning rule is called {\em symmetric} if 
    for all $a,b$ and $S$ as above it holds that 
    $p_{a,b}(S) = p_{a,b}(D \setminus (S \cup \{a,b\}))$. Equivalently, it is symmetric if $p_{a,b}(S) + p_{b,a}(S) = 1$.
\end{defn}

\begin{example}\label{example: wr1}
    Let $D = [2n] = \{1,\ldots,2n\}$ and $p_{a,b}(S) = \frac{1}{2}$ for every $a, b$ and $S$, then on every turn each player has a $\frac{1}{2}$ chance of winning no matter which cards are played. This symmetric winning rule encodes a simple random walk. 
\end{example}
\begin{example}\label{example: wr2}
    Let $D = [2n] = \{1,\ldots,2n\}$, then the winning rule $\mfp$ with $p_{a, b}(S) = \ind\{a > b\}$, the indicator that $a > b$, recovers the simple variant of war mentioned in the introduction. This is a symmetric rule that says in any round, the player with the larger card wins with probability 1. 
\end{example}
Oftentimes, as in the above examples, $\mfp$ does not depend on $S$, the rest of the hand of player's $A$ and $B$. The following is an example which demonstrates the advantage of allowing such a dependence. 

\begin{example}\label{example: wr3}
    Let $D$ be a deck of $2n$ not necessarily distinct cards, with the symmetric winning rule 
    \[ p_{a,b}(S) = \frac{a^s}{a^s + b^s}, \] 
    where $s = \min(|S|,|D \setminus (S\cup \{a,b\})|)$. Here the chance that the higher card wins increases as $s$ gets larger. For example, it is high when the players have about the same number of cards but it is $\frac{1}{2}$ if one of the players is playing their last card. 
\end{example}

We now give an example of a non-symmetric rule to illustrate the importance of symmetry in the behavior of the game. 
\begin{example}\label{example: wr4}
    Let $D$ be a deck of $2n$ distinct cards, with the winning rule 
    $p_{a,b}(S) = 1$ if 
    $\max(S \cup \{a\}) > \max(D \setminus (S \cup \{ a \}))$ and 0 otherwise. 
    This winning rule is not symmetric, and it is easy to see that the player with the highest card wins the game in at most $2 n$ rounds. This differs greatly from the behavior of symmetric winning rules as shown in Theorem \ref{t:deltawar}.
\end{example}    
Some some further examples of winning rules are listed in Corollary~\ref{cor: variants}. The following is a formal description of the game $\mfp$-war from the beginning of this section using the language of Markov chains (see \cite{levin2017markov} for a comprehensive introduction to Markov chains). 

\begin{model}[$\mfp$-war]
  For a fixed a winning rule $\mfp(S) = (p_{a, b}(S))_{a, b \in D}$, $\mfp$-war is a Markov chain on the state space $\Omega = \{(A,B):A,B\subseteq D, A\sqcup B=D\}$, where a state represents the cards which the players $A$ and $B$ have. 
We use $A(t), B(t)$ to denote the hands of the players $A$ and $B$ at time $t$.

\textit{Transition probabilities:} The evolution of $\mfp$-war is given by the following transition probabilities. The absorbing states of this chain are $(\emptyset,D), (D,\emptyset)$ corresponding the player $B$ winning and player $A$ winning. If we are not at an absorbing state, we first draw independent uniformly random cards $a$ in $A(t)$ and  $b$ in $B(t)$ from the hands at time $t \in \NN$ and subsequently sample  $U$ which is a uniform real number between $0$ and $1$ (whose distribution is denoted by $\text{Unif}([0, 1])$).
Then player $A$ wins both cards with probability $p_{a, b}$ and player $B$ wins both cards with the remaining probability $1-p_{a, b}$. 
More formally: 
$$
(A(t+1),B(t+1)) = \begin{cases}
(A(t) \cup \{b\}, B(t) \setminus \{b\}) & U \le p_{a, b}(A(t) \setminus\{a\}), \\
(A(t) \setminus \{a\}, B(t) \cup \{a\}) & \text{otherwise}. \\
\end{cases}
$$

We use $\tau_\mathrm{abs}(\mathfrak{p})$ to denote the time at which the chain absorbs, i.e., one of the players wins the game. For $t \in \NN$, let $A_t$ be the number of cards that player $A$ has at time $t$ and $B_t$ be the number of cards that player $B$ has at time $t$.
\end{model}

Our main results established as long as cards are properly shuffled at the beginning of the game, the number of cards in player $A$'s hand is exactly a simple random walk: 

\begin{thm}\label{t:deltawar}
Consider a game of $\mfp$-war with a deck of $2n$ cards with ranks in $[2n]$ (possibly including repeated ranks) where $\mfp$ is a symmetric winning rule. Suppose that $A(0)$ is uniformly distributed among all hands of size $A_0$. Then, $A_t$ is a simple symmetric random walk stopped at 0 and $2n$.
\end{thm}
\begin{proof}
Consider a game of $\mfp$-war. Suppose that $A(t)$ is uniformly distributed across all possible hands of size $k$. 

Let $a$ and $b$ be the cards played by $A$ and $B$ in the next round, and define $S=A(t)\setminus \{a\}$. By the assumption of $A(t)$ being uniformly distributed across all hands of size $k$, it follows that $S$ is uniform across all hands of size $k-1$. 
Note that once the round is over, both $a,b$ are given to player $A$ with probability $1/2$ or instead to player $B$ with probability $1/2$.

We now describe the above procedure in an alternate but equivalent way. We first reveal the multi-set $S$, which is a uniform hand of size $k-1$. Having done so, we reveal the set $\{a,b\}$ which are a uniform pair of cards chosen from the remaining $2n-k+1$ cards. Finally, we reveal if $A$ won or $B$ won and both are equally likely since only the set $\{a,b\}$ has been revealed, but the individual cards $a,b$ have not; indeed, this step uses the symmetry that $(p_{a,b}(S)+p_{b,a}(S))/2=1/2$ from the definition of a winning rule.
Each valid outcome of the above three step procedure has the probability
\[ \frac{1}{\binom{2n}{k-1}} \cdot \frac{1}{\binom{2n-k+1}{2}} \cdot \frac{1}{2} = \frac{1}{\binom{2n}{k}} \cdot \frac{1}{k} \cdot \frac{1}{2n-k}, \]
and as one would expect, the right hand side agrees with the probability of choosing the hand $A(t)$ uniformly, and then selecting a single card from each of $A(t)$ and $B(t)$ uniformly. 
As a consequence, with equal probability, $A(t+1)$ has either $k-1$ or $k+1$ cards. Moreover, $S$ and $S \cup \{a,b\}$ are uniformly random hands of size $k-1$ and $k+1$ respectively, so $A(t+1)$ remains uniformly random among hands of the given size. 
This implies that with $\frac12$ probability, $A(t+1)$ is uniformly distributed on hands of size $k+1$ and by symmetry, with $\frac12$ probability, $A(t+1)$ is uniformly distributed on the possible hands of size $k-1$. Inductively the desired conclusion follows.
\end{proof}

As a consequence of the above theorem, we can transfer known results for simple random walks (see e.g. \cite[Chapter 6]{van2016probability}) to our setting of $\mfp$-war. 
\begin{cor}\label{cor: variants}
Consider a game of $\mfp$-war with a deck of $2n$ cards with ranks in $[2n]$ (possibly including repeated ranks) where $\mfp$ is a symmetric winning rule. Then: 
\begin{itemize}
\item The expected time for termination is given by: 
$\EE[\tau_{\ab}(\mfp)]  = A_0 B_0$
\item The probability that player $A$
 wins is: $\P(A_{\tau_{\ab}(\mfp)} = 2n) = \frac{A_0}{2n}$. 
\end{itemize}
\end{cor}
In particular, the corollary applies to the following variants of war: 
\begin{itemize}
    \item From Example \ref{example: wr1}, a game of war with $2n$ distinct cards (winning rule $p_{a,b} = \mathbf{1}\{a > b\}$ independent of $S$)
    \item From Example \ref{example: wr2}, a game of war with $2n$ not necessarily distinct cards, where ties are broken by coin flip (winning rule $p_{a,b} = \mathbf{1}\{a > b\} + \frac12 \cdot \mathbf{1}\{a=b\}$ independent of $S$).
    \item From Example \ref{example: wr3}, a game of war with $2n$ not necessarily distinct cards, with the winning rule $p_{a,b}(S) = \frac{a^s}{a^s + b^s}$, 
    where $s = \min(|S|,|D \setminus (S\cup \{a,b\})|)$.
    \item A \emph{martingale} war, where cards are played randomly from each player's hand (winning rule $p_{a,b} = \frac{a}{a+b}$ independent of $S$). The inspiration for the name martingale war will be discussed further in the next section.
\end{itemize}

\section{Deterministic gameplay}\label{s:gladiator}
In the above models of war, we assumed that cards are played randomly from a player's hand. However, in \textit{classic war}, players initially receive a shuffled half of the deck, but thereafter only play the top card(s) from their hand and return cards they win to the bottom of the hand.
Importantly, the manner in which cards are returned to the bottom of the hand can have enormous worst-case implications for the time to win. In fact, it is possible for a game of war with either a traditional 52-card deck or a deck of cards indexed by $[n]$ to have infinite loops if players play their top card and are allowed to return cards to the bottom in arbitrary order (as observed in~\cite{LR12}).

To study the impact of card ordering, in this section, we consider a different variant of war, where we introduce randomness in the outcome of each round, but play cards from the top of each player's hand, and return cards to the bottom of the winning player's hand. The randomness is introduced by selecting a winning card in each round according to a specified Bradley-Terry Model of interaction. This gameplay of the variant proceeds as follows. In each round, player $A$ plays their \textit{top} card $a$ and player $B$ plays their \textit{top} card $b$. Each of these cards has a pre-determined strength $f(a)$ and $f(b)$ respectively. Each player then wins the round with probability proportional to the strength of their respectively played card. The winning player receives both cards, and places them at the bottom of their hand in random order. More formally we have the following description.

\begin{model}[Martingale war]
Fix a function $f: [n] \rightarrow \RR_+$, which is to be thought of as the strength of each card. Let \textit{$f$-war} be a Markov chain defined on state space $\Omega = \{(A, B) : A \sqcup B = [n]\}$, viewing $(A, B) \in \Omega$ as a pair of disjoint \textit{ordered} tuples of elements of $[n]$.

Any state with $|A| = n$ or $|B| = n$ is an absorbing state, i.e. a situation where either player $A$ or $B$ wins. We let $(A(t), B(t))$ denote the hands of players $A$ and $B$ at time $t \in \NN$. 
We now describe the transition probabilities. If $(A(t), B(t))$ is not an absorbing state, then $(A(t), B(t)) \mapsto (A({t+1}), B({t+1}))$ as follows: let $a$ be the first card in $A(t)$ and $b$ be the first card in $B(t)$. Then, with probability $\frac{f(a)}{f(a) + f(b)}$, we add $a$ and $b$ in arbitrary order to the bottom of $A(t) \backslash \{a\}$ to construct $A(t+1)$ and let $B(t+1) = B(t) \backslash \{b\}$. Otherwise, we add $a$ and $b$ in arbitrary order to the bottom of $B(t) \backslash \{b\}$ to construct $B(t+1)$ and let $A(t+1) = A(t) \backslash \{a\}$.
\end{model}

To motivate the naming of this variant, we introduce the martingale and the properties relevant to our analysis. 

\subsection{Martingales}
A martingale is a probabilistic object that represents the process of playing a fair game, meaning that the average winnings from playing the game are zero. Ville \cite{ville1939etude} gave this name to a concept introduced by L\'evy to study sequences of dependent random variables, following a rich history of using this term informally to describe betting strategies \cite{mansuy2009origins}. Martingales have since been used to study betting strategies in games of chance and are prevalent in probability. A short introduction (in the discrete setting) is below, with a simple example to follow.

\begin{defn}[Martingale]
A \textit{martingale} is a sequence of random variables $\{X_t\}_{t \in \NN}$ such that for all $t \in \NN$, $X_t$ is finite on average and moreover if we know the history up to time $t$, then the average value of $X_{t+1}$ is exactly $X_t$. More formally, this is expressed by: $\EE[|X_t|] < \infty$ and $\EE[X_{t+1} \mid X_0, \ldots, X_t] = X_t$ for all $t$.
\end{defn}

A fundamental result in the theory of martingales is the Optional Stopping Theorem.
\begin{thm}[Optional Stopping Theorem]
    Let $\{X_t\}_{t \in \NN}$ be a martingale such that there exists $K$ with $|X_t - X_{t-1}| < K$ uniformly for all $t$. Let $\tau$ be a (random) time such that $\EE[\tau] < \infty$. Then $\EE[X_\tau] = \EE[X_0]$. 
\end{thm}
The theorem says that after a (potentially random) amount of time, the average value of the martingale is the same as the average initial value. It is very useful in studying properties of these random times. Notice that the definition of a martingale implies $\EE[X_t] = \EE[X_0]$ for all fixed times $t$. The optional stopping theorem gives conditions under which this extends to a certain class of random times $\tau$ called stopping times. Being a stopping time roughly means that the decision to stop at time $\tau=t$ depends only on $t$ and $X_0,\ldots,X_t$. The conditions of the theorem written above are specialized to our application, and are sufficient but not necessary. See e.g. \cite{durrett2019probability} for further details, such as the general conditions for the optional stopping theorem, and a comprehensive introduction to martingales. 

\begin{example}\label{example: rw}
    Consider one of the simplest examples of a martingale. Let $S_t$ the simple symmetric random walk started at 0, where at each time step we add or subtract 1 with equal probability.
    One can check that $S_t^2 - t$ is also a martingale (one such proof is in \cite{durrett2019probability}). Now let $\tau$ be the first time that the random walk reaches either $n$ or $-n$. The Optional Stopping Theorem implies $\EE[S_\tau^2 - \tau] = \EE[S_0 - 0] = 0$. Notice that $S_\tau^2 = n^2$ regardless of whether the random walk arrives at $n$ or $-n$ first. All together, $\EE[\tau] = n^2$, which recovers the quadratic termination time alluded to in the introduction. See e.g. \cite{van2016probability} for the technical details have been hidden in this example, including the conditions of the Optional Stopping Theorem. 
\end{example}

\subsection{Analysis of martingale war}
In the context of martingale war, the sequence of random variables $M_t^{(f)}$ is a martingale, where 
$$M_t^{(f)} := \sum_{a\in A(t)}f(a), \quad t \in \NN$$
is the total strength of player A's hand at time $t$. Below we list several examples of the strength function $f$ and the resulting war games. 

\begin{example}\hfill
\begin{itemize}
    \item When $f(a) = 1$ is a constant function, $f$-war is a simple random walk. 
    \item When $f(a) = a$ is the identity function, $f$-war is the Gladiator game (see \cite{winkler2003mathematical}).
    \item When $f(a) = e^{\lambda a}$, the game play can be tuned by a parameter $\lambda$, and approaches standard war when $\lambda \to \infty$. 
\end{itemize}
\end{example}

\begin{rem}
Applying the Optional Stopping Theorem to $M_t^{(f)}$ has implications on the probability that each player wins. Suppose that at the beginning, card $n$ is given to player $A$ and every other card is in randomly and independently given to player $A$ with probability $1/2$. Let $p_A^n$ be the probability that $A$ wins with this initial distribution. By the Optional Stopping Theorem,
\[ f(n)+\frac12\sum_{i=1}^{n-1} f(i)= \EE[M_0^{(f)}] = \EE[M_{\tau_{\ab}}^{(f)}] = p_A^n\sum_{i=1}^n f(i). \]
Solving for $p_A^n$ gives
\begin{equation}
  \label{eq:5}
  p_A^n=\frac{1}{2}+\frac{f(n)}{2\sum_{i=1}^nf(i)}.
\end{equation}
For example if $f(a) = a$, the probability that player $A$ wins is only $p_A^n = \frac12 + \frac{1}{n+1}$. This behavior is quite different from classic war, where having the card $n$ guarantees winning. However, as noted above, for $f(a) = e^{\lambda a}$, $p_A^n \to 1$ as $\lambda \to \infty$, and the game approaches the usual war, where the higher card wins on each turn.
\end{rem}

As demonstrated in Example \ref{example: rw}, finding a martingale based on $S_t^2$ allows one to compute the expected hitting time of the simple random walk. The next lemma gives the analogous martingale for $(M_t^{(f)})^2$. 
\begin{lemma}\label{lemma: doob}
For $t \in \NN$, let $Q_t^{(f)} = \sum_{k=1}^t f(a_k)f(b_k)$, where $a_t$ and $b_t$ are the top cards of $A(t)$ and $B(t)$, respectively. Then,
$\{(M_t^{(f)})^2 - Q_t^{(f)}\}_{t \in \NN}$ 
is a martingale. 
\end{lemma}
\begin{proof}
We compute the average value of $(M_{k+1}^{(f)})^2$ given the history of the hand strength. 
  \begin{align*}
      \EE\left[\left(M_{k+1}^{(f)}\right)^2 \middle\vert M_{k}^{(f)}\right] &= \frac{f(a_k)}{f(a_k) + f(b_k)}\cdot (M_k^{(f)} + f(b_k))^2 + \frac{f(b_k)}{f(a_k) + f(b_k)}\cdot (M_k^{(f)} - f(a_k))^2 \\
      &= (M_k^{(f)})^2 + \frac{f(a_k)f(b_k)^2}{f(a_k) + f(b_k)} + \frac{f(a_k)^2f(b_k)}{f(a_k) + f(b_k)} \\
      &= (M_k^{(f)})^2 + f(a_k) f(b_k).
  \end{align*}
Observe that
  \[ Q_t^{(f)} = \sum_{k=1}^t f(a_k)f(b_k) = \sum_{k=1}^t \left( \EE[(M_k^{(f)})^2 \vert M_{k-1}^{(f)}] - (M_{k-1}^{(f)})^2 \right), \]
is precisely the cumulative change in expectation up to time $t$. Thus, $\{(M_t^{(f)})^2 - Q_t^{(f)}\}_{t \in \NN}$ 
is a martingale. 
\end{proof}

\begin{rem}
    The above computation is known as finding the \textit{quadratic variation} of the process $(M_t^{(f)})^2$ using the \textit{Doob decomposition}. The quadratic variation and Doob decomposition are central to the theory of stochastic processes, for more on these concepts see e.g. \cite{durrett2019probability} or \cite{le2016brownian}.
\end{rem}

As a consequence of finding the above martingale, we show that a simple variant of war (taking $f(a) = a + n$) with cards being played from the top of each players' hand takes $\Theta(n^2)$ rounds to end on average. 

\begin{claim}
Consider a game of $f$-war with $f(a) = a+n$, and $A(0)$, the initial hand of player $A$, distributed as follows. Let each card in the deck $[n]$ be given to player $A$ with probability $\frac12$. Then randomly permute the cards assigned to player $A$ to obtain $A(0)$.  Let $\tau_\mathrm{abs}$ denote the time until one of the players wins. Then, $ \EE[\tau_\mathrm{abs}] = \Theta(n^2)$. 
\end{claim}
\begin{proof}
  We first compute some preliminary statistics of $M_t$. By the linearity of expectation,
$\EE[M_0] = \sum_{i=n+1}^{2n} \frac{i}{2} = \frac{3n^2+n}{4}.$
We next compute the second moment:
  \begin{align*}
    \EE[M_0^2] &= \EE\left[\left( \sum_{i=n+1}^{2n} i\ind_{i \in A} \right)^2\right] 
    = \frac{1}{4} \sum_{i=n+1}^{2n} i^2 + \frac{1}{4} \sum_{i=n+1}^{2n}\sum_{j=n+1}^{2n} ij 
    = \frac{9n^4}{16} + O(n^3).
  \end{align*}
  By the symmetry between players $A$ and $B$, $\P(A\text{ wins}) = 1/2$. This implies that  
    \[ \EE[M_{\tau_\mathrm{abs}}] = \frac{1}{2} \sum_{i=n+1}^{2n} i + \frac{1}{2} \cdot 0 = \frac{3n^2+n}{4}; \quad \EE[M_\tau^2] = \frac{1}{2} \cdot \left(\frac{3n^2+n}{2}\right)^2 = \frac{9n^4}{8} + O(n^3). \]
  Recall the process $Q_t$ from Lemma \ref{lemma: doob}. Now we would like to apply the optional stopping theorem to $\tau_\mathrm{abs}$. Since $M_t^{(f)}$ is bounded by the total strength of all cards, and the increment of $Q_t$ is the product of the strengths of two cards, we have 
  \[ |(M_t^2 - Q_t) - (M_{t-1}^2 - Q_{t-1})| < 2\left(\sum_{a=1}^n f(a)\right)^2 + \max_{1 \leq a \leq n} f(a)^2, \]
  so the martingale differences are bounded by $K$
  given by the right hand side above. 
  To bound the expected value of $\tau_\mathrm{abs}$, note that in each round the probability 
  that player A wins the round is at least
  $\frac{f(1)}{f(1)+f(n)} \geq \frac{1}{3}$. 
  Thus, in $n$ consecutive rounds, there is a probability of at least $3^{-n}$ that player A wins all rounds in a row until they win the game. It is easy to see and well known that this implies that 
  $\EE[\tau_\mathrm{abs}] \leq 3^n < \infty$. Thus, combining the optional stopping theorem with the above calculations gives
  \begin{align*}
      \EE[Q_{\tau_\mathrm{abs}}] &= \EE[M_{\tau_\mathrm{abs}}^2] - \EE[M_0^2] = \frac{9n^4}{16} + O(n^3).
  \end{align*}
  The increment in $Q_t$ at each time $t$ is given by the values of the current cards in play. The contribution at each step is at least $n^2$, so 
  \[ \EE\left[ n^2 \tau_\mathrm{abs} \right] \leq \EE[Q_{\tau_\mathrm{abs}}] = \frac{9n^4}{16} + O(n^3),\]
  implying the bound $\EE[\tau_\mathrm{abs}] = O(n^2).$
  Conversely, the contribution at each step is at most $4n^2$, so 
  \[ \EE\left[ 4n^2 \tau_\mathrm{abs} \right] \geq \EE[Q_{\tau_\mathrm{abs}}] = \frac{9n^4}{16} + O(n^3),\]
  which yields the desired lower bound. 
\end{proof}

\begin{rem}
    The same proof applied to the gladiator game with $f(a) = a$ yields an upper bound of $O(n^3)$, whereas we believe the true upper bound should be $O(n^2)$. Similarly, for functions $f(a) = e^{\lambda a}$, the obtained exponential upper bound should not be tight either.
\end{rem}

\section{Simulations and Open questions} \label{s:variants}

\subsection{Simulations with a standard deck of cards} \label{s:experiment}
We compare the above bounds to empirical results obtained by simulating the game of war with a standard $52$-card deck, breaking numerical ties via a \textit{war round}, 
where players flip over their next card (ignoring its value), using their next next card from the top to break the tie.
This rule is iterated as needed. 
A player loses if they run out of cards in the middle of breaking the tie. 

\begin{figure}[ht!]
    \captionsetup{justification=centering,margin=2cm}
    \centering
    \begin{subfigure}[b]{0.47\textwidth}
    \centering
    \includegraphics[width=0.97\linewidth]{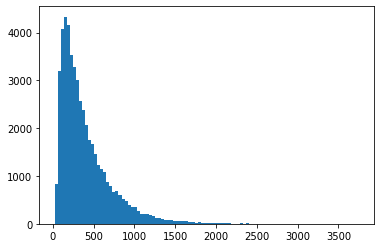}
    \caption{We simulate war with a standard $52$-card, 4-suit deck, drawing cards from the top of each player's hand and breaking ties via a \textit{war round}. %(see above).
    \\Mean \#rounds: $397$ \\Median \#rounds: $302$ \\Max \#rounds: $3752$}
    \end{subfigure}
    \begin{subfigure}[b]{0.47\textwidth}
    \centering
    \includegraphics[width=0.97\linewidth]{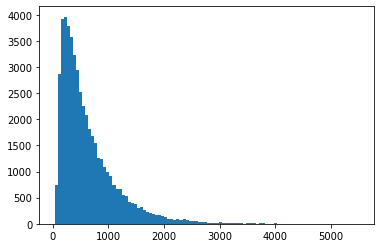}
    \caption{We simulate war with a standard $52$-card, 4-suit deck, drawing cards from the top of each player's hand and breaking ties by a fair coin flip. \\Mean \#rounds: 628 \\Median \#rounds: 474 \\ Max \#rounds: 5510}
    \end{subfigure}
    \begin{subfigure}[b]{0.47\textwidth}
    \centering
    \includegraphics[width=0.97\linewidth]{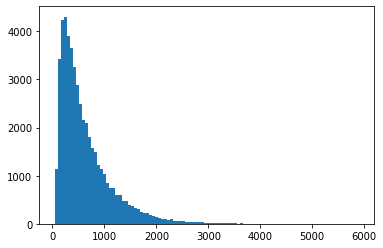}
    \caption{We simulate war with a standard $52$-card, 4-suit deck, drawing cards \textit{randomly} from each player's hand and breaking ties by a fair coin flip. \\
    Mean \#rounds: 625 \\
    Median \#rounds: 472 \\
    Max \#rounds: 5900}
    \label{fig:histmodel}
    \end{subfigure}
    \begin{subfigure}[b]{0.47\textwidth}
    \centering
    \includegraphics[width=0.97\linewidth]{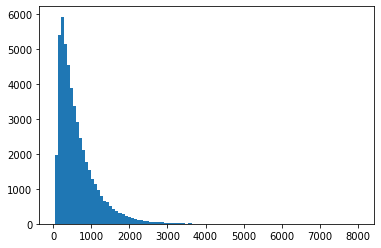}
    \caption{We simulate war with a deck of $52$ cards, each of distinct rank, drawing cards from the top of each player's hand. \\
    Mean \#rounds: 624 \\
    Median \#rounds: 474 \\
    Max \#rounds: 8026}
    \end{subfigure}
    \caption{Histograms of the number of rounds in three models of war}
    \label{fig:hists}
\end{figure}

We empirically compare four different models of war with a $52$-card deck;~\cref{fig:hists} depicts four histograms of the number of rounds in the associated games (we simulated 50,000 games for each model). Earlier, we bounded the expected winning times of the model pictured in~\cref{fig:histmodel}. 

When playing war with a single-suited deck of cards, each of distinct rank, the person with the strongest card is \textit{guaranteed} to win. This is not the case with a standard deck of cards (or any deck with repeated cards of the strongest rank). We experimentally compute (with a standard deck of cards) the probability of winning a game of war conditioned on the number of cards of the strongest rank a player starts with (which we assume is an Ace). These probabilities are given below; in both models, we use a standard $52$-card, $4$-suited deck and have each player draw cards from the top of their hands, but in model (A), we break ties via a war round, and in model (B), we break ties by flipping a fair coin.

\begin{center}
\begin{tabular}{|c|c|c|} 
\hline 
  \# of Aces   & $\P$(win) in model (A) & $\P$(win) in model (B) \\
  \hline 
  \hline 
  0   &  0.108 & 0.000 \\
  1   &  0.293 & 0.243 \\
  2   &  0.500 & 0.500\\
  3   &  0.706 & 0.757\\
  4   &  0.892 & 1.000\\
  \hline 
\end{tabular}
\end{center}

The above simulation highlights that \textit{war rounds} cause a player with no Aces to still have a relatively large ($\approx 10\%$) chance of winning a game of war (as seen in model (A)); this is in stark in contrast to model (B), where a player can \textit{never} win unless they start with at least one Ace.

\subsection{Open questions}
Of course, a primary question left open by this work is understanding the expected winning time of a game of \textit{classic war}, played with a standard deck of $mn$ cards with $n$ distinct ranks and $m$ copies of each rank (for constant $m$), where (a) players draw cards from the top of their hands each round, as opposed to randomly and (b) ties are broken by \textit{war rounds}, rather than coin flips. From simulations (as summarized by~\cref{fig:hists}), classic war appears to terminate much more quickly than \textit{random war} the variant we primarily study in this work. It appears heuristically as though the \textit{war} rounds, which cause many additional cards to change hands in a single round, speeds up play in the classic game. A calculation reveals that this speed up is quantified by an increase in quadratic variation (for example when a random walk can take steps of size 3).

Nonetheless, we conjecture that the expected winning time remains $O(n^2)$ (similarly to random war) as the scaling effects seem similar when simulating both games with varying hand size.
Exhibiting $O(n^2)$ expected winning time under assumption (a) would already be of interest.
It is similarly interesting to extend the result presented here to cases where there are more than two players. 

Much of this work focused on studying the expected time to win in a game of \textit{random war}. A different, natural question is to understand the \textit{probability} that a player wins a game of war, given some partial information about their hand composition; for example, the above numerical simulations estimate the probability of winning a game of war given the number of Aces in a player's starting hand. Understanding both what types of partial information are most predictive of the chance a player wins a game of war (e.g. the sum of cards, the number of cards of some specific ranks, $\ldots$ etc) and the conditional probability of winning given some summary statistics of a player's hand are both interesting problems.

\subsection*{Acknowledgements}
MB was supported by a Peterson fellowship at MIT. BC was supported by the NSF Graduate Research Fellowship Program. NM was supported by a Hertz Graduate Fellowship and the NSF Graduate Research Fellowship Program. The questions raised in this paper were inspired by many war games with Nathan Mossel. E.M. thanks Nathan for the inspiration in particular by insisting to play with all Aces initially assigned to him. E.M. 
is partially supported by Vannevar Bush Faculty Fellowship award ONR-N00014-20-1-2826 and by a Simons Investigator Award in Mathematics (622132).

\bibliographystyle{amsplain0}
\bibliography{errorcorrecting.bib}

\end{document}